\documentclass[12pt]{article}

\usepackage{amssymb,amsmath,amsfonts} 
\newtheorem{theorem}{Theorem} 
 
\newtheorem{proposition}{Proposition} 
\newtheorem{corollary}{Corollary} 
\newtheorem{definition}{Definition}

\begin{document}

\title{ On the Lifts of Minimal Lagrangian Submanifolds } 

\author{ Sung Ho Wang  }
\date{September, 2001}
\maketitle

Bryant and Salamon constructed metrics with holonomy $G_2$ and $Spin(7)$ on 
spin bundles of 3-dimensional space forms, and spin bundles and bundles of anti-self-dual
2-forms on self-dual Einstein 4-manifolds [BrS]. Since, apart from holonomy, the construction of 
integrable $G_2$(respectively $Spin(7)$) structures amounts to finding differential 3(4)-forms of 
generic type on $7$($8$) manifolds satisfying appropriate differential equations, the tautological and
canonical forms on these bundles serve as the basis of their construction. 

In this paper, we show the total space of the canonical line bundle $\mathbb{L}$ of a Kahler-Einstein manifold $X^n$ supports integrable
$SU(n+1)$ structures, or Calabi-Yau structures. The parallel holomorphic volume form in each case is the exterior derivative
of the canonical $(n, 0)$ form on $\mathbb{L}$.

Since each of integrable $SU(n+1)$, $G_2$, and $Spin(7)$ structures gives rise to a calibration on the underlying manifold,
a question arises as to what the calibrated submanifolds are with regard to the above constructions. We show in $SU(n+1)$ case,
the canonical real line bundle $L \subset \mathbb{L}$ over a minimal Lagrangian submanifold $M \subset X$ is calibrated and
hence can be considered as the special Lagrangian lift of $M$. As a corollary, we show the area of compact connected minimal
Lagrangian submanifolds in a complex projective space with the standard Kahler structure admits a uniform lower bound. 

In $G_2$, and $Spin(7)$ cases, minimal surfaces in self-dual Einstein 4-manifolds with 
vanishing complex quartic form (super-minimal) admit lifts which are calibrated, i.e., associative, coassociative or Cayley respectively. The lifts in 
this case can be considered as the tangential lifts or normal lifts of the minimal surface adapted to the quaternionic bundle structure.  
Compact minimal surfaces with vanishing complex quartic form are known to exist in $S^4$ and $\mathbb{C}P^2$ in abundance from the work of Bryant, 
and Chern and Wolfson [Br1][ChW].

\section{ Cannonical Line Bundle }

Let $X$ be a Kahler manifold of dimension $n$ with metric $g_0$, Kahler form
$\varpi$ such that 
\begin{equation}
g_0(x,y) = \varpi(x, Jy)
\end{equation}
where $J$ is the complex structure on $X$ and $x, y \in TX$. Let $\pi : F \to X$ be the associated $U(n)$
bundle with tautological forms $\eta^k = \omega^k + i \theta^k, k=1, .. n$ for which
\begin{equation}
\pi^* (\varpi) = \omega^1 \wedge \theta^1 + \, . . . \, \omega^n \wedge \theta^n.
\end{equation}
There exists a unique set of connection 1-forms $\alpha^j_k = - \alpha^k_j, \beta^j_k = \beta^k_j$
on $F$ so that the following structure equation hold.
\begin{align}
d\omega^j &= - \alpha^j_k \wedge \omega^k + \beta ^j_k \wedge \theta^k \notag \\
d\theta^j &= - \beta^j_k \wedge \omega^k - \alpha ^j_k  \wedge \theta^k \notag.
\end{align}
We denote
\begin{equation}
\gamma = \textnormal{tr} \beta = \sum_k \beta^k_k
\end{equation}
and note the identity 
\begin{equation}
d\gamma = \pi^* (Ric)
\end{equation}
where $Ric$ is the Ricci form. $X$ is called Kahler-Einstein if
\begin{equation}
 Ric = c \, \varpi,
\end{equation}
for some function $c$, which is necessarily constant when $n \geq 2$. 

The canonical line bundle
\begin{equation}
\mbox{$\bigwedge^{n, 0}$} T^*X = \mathbb{L} \to X
\end{equation}
is the bundle of $(n, 0)$ forms on $X$. Let $\Upsilon_0$ be the tautological form on $\mathbb{L}$ and 
\begin{equation}
r : \, \mathbb{L} \to [0, \infty)  
\end{equation}
denote the norm $r(u) = |u|$ on $\mathbb{L}$. We denote $B = r^{-1}(1) \subset \mathbb{L}$ 
the associated unit circle bundle. By definition, $B = F/SU(n)$ and the tautological $(n,0)$ form $\Upsilon_0$ on $\mathbb{L}$ 
satisfies the following structure equation when restricted to $B$.
\begin{align}
\pi^*(\Upsilon_0|_{B})&=\frac{1}{2^{\frac{n}{2}}} \eta^1 \wedge \,  . . . \, \wedge \eta^n  \\
d \Upsilon_0 &= - i \, \gamma \wedge \Upsilon_0  \\
d \gamma &= Ric = c \, \varpi.
\end{align}
where $\pi : F \to B$.

\section{Calabi - Yau Structure}
Set
\begin{equation}
\Upsilon = d \Upsilon_0.
\end{equation}
Then $\Upsilon$ is a closed $(n+1)$ form on $\mathbb{L}$ which is holomorphic with respect to the
underlying complex structure. We wish to show there exist Kahler structures on the total space
of the line bundle $\mathbb{L} \to X$ under which $\Upsilon$ is a holomorphic volume form
of constant length (hence parallel). For simplicity, we assume $n \geq 2$.

Consider a parametrization
\begin{align}
\varphi: R^+ \times B &\to \mathbb{L} \\
(r, u) &\to r \, u
\end{align}
and set
\begin{equation}
\eta^0 = dr - i \,  r \, \gamma.
\end{equation}
Then we have
\begin{equation}
\varphi^*(\Upsilon)  = \eta^0 \wedge \Upsilon_0
\end{equation}
and $\{ \eta^0, \eta^1, \, . . \, \eta^n \}$ is  a basis of $(1,0)$ forms on $\mathbb{L}$.

The metric $g$ on $\mathbb{L}$ we are interested in is of the form
\begin{align}
g &= f^{-2}(r) \eta^0 \cdot \overline{\eta^0} + f^{\frac{2}{n}}(r) 
            (\eta^1 \cdot \overline{\eta^1} + \, . . . \, \eta^n \cdot \overline{\eta^n}) \\
\Pi &= \frac{i}{2} ( f^{-2}(r) \eta^0 \wedge \overline{\eta^0} + f^{\frac{2}{n}}(r) \varpi )
\end{align}
where $\Pi$ is the associated 2-form and $f(r) > 0$. Note that under this Hermitian structure, $\Upsilon$ is
a closed $(n+1,0)$ form of constant length and hence it will be parallel whenever the Hermitian
structure is Kahler, i.e., $\Pi$ is closed. 

By (10) and (14), $d\Pi = 0$ is equivalent to
\begin{equation}
c \, r + \frac{2}{n} f^{\frac{2}{n}+1} \, \frac{\partial{f}}{\partial{r}} = 0
\end{equation}
or
\begin{equation}
f(r) = ( -c \frac{(n+2)}{2} r^2 + c' )^{\frac{n}{2n+2}}
\end{equation}
where $c'$ is an arbitrary constant.

\begin{theorem}
Let $\mathbb{L} \to X$ be the canonical line bundle of a Kahler-Einstein manifold
$X$. Then the metric (16) with $f(r)$ in (19)  gives rise to a Calabi-Yau structure on
an open set of $\mathbb{L}$. If $c \leq 0$ and $c' >0$, the metric is defined on all of 
$\mathbb{L}$ and it is complete if the original metric on $X$ is complete.
\end{theorem}

Note that these Calabi-Yau structures provide calibrating forms  $2^{\frac{n}{2}} Re(e^{i \theta} \Upsilon)$.

\section{Lifts of Minimal Lagrangian Submanifolds}

Let $\phi: M \to X$ be a Lagrangian submanifold of a Kahler-Einstein manifold $X$. From (10),
the $S^1$ bundle $\phi^*(B) \to M$ is flat. In fact, let $\{ \, e_1, \, . . . \, e_n \, \}$ be 
a local basis of $TM$ and set $n_k = J e_k$, and consider 
\begin{equation}
s = \frac{1}{2^{\frac{n}{2}}} (e_1^* + i n_1^*) \wedge \,  . . . \, \wedge (e_n^* + i n_n^*)
\end{equation}
where $*$ denotes the metric dual. Then  $s$ is a well defined (mod $\mathbb{Z}_2$ if $M$ is not orientable) 
section of $\phi^*(B)$. 
It is well known, [Br2], that $s$ is parallel, or $s^*(\gamma)=0$, if and only if $M$ is minimal.  

\begin{definition}
Let $\phi: M \to X$ be a connected Lagrangian submanifold of a Kahler manifold $X$. 
Let $B \to X$ be the canonical $S^1$ bundle. A section (or mod $\mathbb{Z}_2$ section if $M$ is not 
orientable) $\, s: M \to B$ is called a $Legendrian \, \, lift$ of $M$ if 
$s^*(\gamma) = 0$, where $\gamma$ is the connection form
of the bundle $B \to X$. We denote the Legendrian lift by $\tilde{M}$.
\end{definition}
Note that $\tilde{M}$ is connected.

\begin{proposition}
Let $M \subset X$ be a connected minimal Lagrangian submanifold of a Kahler-Einstein manifold.
Let $B \to X$ be the canonical $S^1$ bundle. There exists a $S^1$ family of Legendrian lifts $\tilde{M} \subset B$
such that $\tilde{M} \to M$ is 1:1 if $M$ is orientable and 2:1 if $M$ is nonorientable.
\end{proposition}

\begin{theorem}
Let $\phi : M \to X$ be a minimal Lagrangian submanifold of a Kahler-Einstein manifold $X$.
Let $\mathbb{L} \to X$ be the canonical line bundle with tautological form
$\Upsilon_0$ and  Calabi-Yau structure (16).
Let $L \subset \mathbb{L} \to M$ be the smooth real line bundle generated by a Legendrian lift 
$\tilde{M} \subset B \subset \mathbb{L}$ of $M$, which is well defined even if $M$ is not orientable.
Then $L$ is a special Lagrangian submanifold of $\mathbb{L}$ calibrated by $2^{\frac{n}{2}} \, Re(e^{i \theta} d\Upsilon_0)$
for some $e^{i \theta}$.
\end{theorem}
Remark : Consider the zero section of the bundle $\phi^*(\mathbb{L}^{n+1}) \to M^n$, which is an isotropic
submanifold of $\mathbb{L}$. The theorem above provides an $S^1$ pencil of minimal Lagrangian submanifolds
containing the given isotropic submanifold, as observed by Bryant [Br2]. Also note that $L$ is in general 
not the fixed locus of an anti-holomorphic involution.

Proof. $\; $ Up to standard $S^1 \subset \mathbb{C}^*$ action on $\mathbb{L}$, we may assume $\tilde{M}$ 
is given by (20). Now it easily follows from the description of 
Calabi-Yau structure that $L \subset \mathbb{L}$ is a Lagrangian submanifold on which 
$Im(d\Upsilon_0) = 0$. $\square$
\vspace{1pc}

Consider $\mathbb{C}P^n = SU(n+1)/U(n)$, which is  Kahler-Einstein with $c=2(n+1)$. Let $S^{2n+1} \to \mathbb{C}P^n$
be the Hopf map. Then the associated unit $S^1$ bundle in this case is $B = S^{2n+1}/ \mathbb{Z}_{n+1} \to \mathbb{C}P^n$, where 
$\mathbb{Z}_{n+1}$ represent the center
of $SU(n+1)$. Given a connected minimal Lagrangian submanifold $M \subset \mathbb{C}P^n$, it is well known there is a maxiaml connected
lift $\tilde{M} \subset S^{2n+1}$, unique up to standard $S^1$ action, that is minimal and Legendrian.
In view of Theorem 1., we have the following corollary.

\begin{corollary}
Let $\phi : M \to \mathbb{C}P^n$ be a connected minimal Lagrangian submanifold, and let 
\begin{equation}
\pi_1(M) \to H_1(M, \mathbb{Z}) \to Hol(M, S^1) \subset O(2)
\end{equation}
denote the holonomy of the associated flat $S^1$ bundle $\phi^*(S^{2n+1}) \to M$. Then
\begin{align}
Hol(M, S^1) &\subset \mathbb{Z}_{n+1}  \subset SO(2) \textnormal{\; if and only if $M$ is orientable} \\
Hol(M, S^1) &\subset \mathbb{D}_{n+1}  \subset O(2) \textnormal{\; if $M$ is nonorientable} 
\end{align}
where $\mathbb{D}_{n+1}$ is the difedral group of order $2(n+1)$.
In particular, if $M$ is compact or embedded, there exists a connected minimal Legendrian lift $\tilde{M} \subset S^{2n+1}$
that is compact or embedded respectively.
\end{corollary}

Proof. $\; $ Take a point $p \in M$ and let $RP^n_p$ be the unique tangent plane to $M$ at $p$. 
It is well known that if $\tilde{M} \subset S^{2n+1}$ is a minimal Legendrian lift of $M$, then
the cone $0 \times \tilde{M}$ is minimal Lagrangian in $\mathbb{C}^{n+1}$. Since $\tilde{M}$ is connected, we 
may assume 
$0 \times \tilde{M}$ is special Lagrangian in $\mathbb{C}^{n+1}$. Consider the special Lagrangian
$(n+1)$ planes in $\mathbb{C}^{n+1}$ which is mapped to $RP^N_p$ under the projection $\mathbb{C}^{n+1}-\{0\} \to \mathbb{C}P^n$.
It is easy to see there  are $2 (n+1)$ of them, $(n+1)$ if we ignore the orientation. Since the inverse image
of a point of the Hopf map is $S^1$, this implies given $p \in M$, there is at most $2(n+1)$ 
inverse images of $p$ in $\tilde{M}$. $\square$

\vspace{1pc}

\textbf{Example 1}
Consider $0 \times SU(n) \subset M(n, \mathbb{C}) = \mathbb{C}^{n^2}$, which is a  minimal Lagrangian cone. Hence
$SU(n)/\mathbb{Z}_n \to \mathbb{C}P^{n^2-1}$ is a minimal Lagrangian embedding. The induced $S^1$ bundle
has holonomy $\mathbb{Z}_n \subset \mathbb{Z}_{n^2}$.

\textbf{Example 2}
  Consider $0 \times Slag(n) \subset Sym(n, \mathbb{C}) = \mathbb{C}^{n(n+1)/2}$, where 
\begin{align}
Slag(n) &= SU(n)/SO(n) \\
        &= \{ \, A A^t \in Sym(n, \mathbb{C} ) \, | \; A \in SU(n) \; \}.
\end{align}
It is a minimal Lagrangian cone and hence its image under Hopf map is a 
minimal Lagrangian submanifold of $\mathbb{C}P^{n(n+1)/2-1}$.
The holonomy of the induced $S^1$ bundle is  $\mathbb{Z}_n \subset \mathbb{Z}_{n(n+1)/2}$ if 
$n$ is odd and  $\mathbb{Z}_{\frac{n}{2}} \subset \mathbb{Z}_{n(n+1)/2}$ if $n$ is even.

\textbf{Example 3}
   Consider the Hexagonal torus
\begin{equation}
T^n = \{ \, [ z^0, z^1, \, . . . \, , z^n ] \in \mathbb{C}P^n \, | \, |z^i| = |z^j| \; \textnormal{ for all \, } i, j   \, \}
\end{equation}
It has a Legendrian lift in $S^{2n+1}$ which is again a torus and an $(n+1)$ fold cover, i.e., $Hol(T^n, S^1) = \mathbb{Z}_{n+1}$.

\begin{corollary}
Let $M$ be a compact minimal Lagrangian submanifold of $\mathbb{C}P^n$. There exists a positive constant $\delta(n)$ such that
Area($M$) $\geq \delta(n)$.
\end{corollary}

Proof. $ \;$ Consider a compact minimal Legendrian submanifold $\tilde{M} \subset S^{2n+1}$. Since the cone over
$\tilde{M}$ is absolutely minimizing in $\mathbb{C}^{n+1}$, isoperimetric inequality provides a uniform area lower bound for such $\tilde{M}$. 
The corollary follows from Corollary 1., for Hopf map is an isometric submersion. $\square$

\section{ Associative and Coassociative Lifts}

Let $N$ be a self-dual Einstein 4-manifold and let $\pi : \bigwedge^2_{-} \to N$ be the bundle of anti-self-dual
2-forms and $ \, \mathbb{S}$ $\to N$ be the spin bundle of real rank 4. 

Given a surface $ \phi : \Sigma^2 \to N$, consider
\begin{align}
\Omega^{\perp} &= \{ \, \xi \in \mbox{$\bigwedge^2_{-}$} \, |  \; \pi(\xi) = p \in \Sigma, \, \, \xi(T_p \Sigma) = 0 \, \} \\
\Omega &= (\Omega^{\perp})^{\perp}.
\end{align}
Then $\Omega$ and  $\Omega^{\perp} \subset \bigwedge^2_{-}$ are a real line bundel and a rank two vector bundle over $\Sigma$ respectively.
For the spin bundle $\mathbb{S}$, set $\tilde{T}\Sigma$ and $\tilde{N}\Sigma$ be the appropriate lifts in $\,\mathbb{S}$ of 
$T\Sigma$ and $N\Sigma$..

A minimal surface in a Riemannian manifold of dimension $\geq 4$ is called superminimal if the associated complex 
quartic form vanishes [ChW].

Bryant and Salamon constructed integrable $G_2$ structures on $\bigwedge^2_{-}$ and integrable $Spin(7)$ structure on $\, \mathbb{S}$
by similar idea used in Section 2 [BrS]. We record the following observation analogous to Theorem 2. without proof.

\begin{theorem}
Let $\Sigma$ be a minimal surface in a self-dual Einstein 4-manifold $N$ with vanishing complex quartic form. Then, under the $G_2$ structures on 
$\bigwedge^2_{-} \to N$  and $Spin(7)$ structures on the spin bundle \textnormal{$\mathbb{S}$} $\to N$ 
constructed by Bryant and Salamon,

1. $\Omega$ is an associative and $\Omega^{\perp}$ is a coassociative submanifold of $\bigwedge^2_{-}$.

2. $\tilde{T}\Sigma$  and  $\, \tilde{N}\Sigma$ are Cayley submanifolds of \textnormal{$\mathbb{S}$}.
\end{theorem}

From the work of Bryant, and Chern and Wolfson, [Br] and [ChW], it is known  $S^4$ and 
$\mathbb{C}P^4$ support many compact minimal surfaces with vanishing
quartic form.

\vspace{5pc}

\noindent
\textbf{\large{References}}

\noindent
[Br1] Bryant, Robert L., Conformal and minimal immersions of compact surfaces
into the $4$-sphere. J. Differential Geom. 17 (1982), no. 3, 455-473

\noindent
[Br2] $\underline{\; \;}$, Minimal Lagrangian submanifolds of Kahler-Einstein
manifolds. Differential geometry and differential equations (Shanghai, 1985), 1--12, Lecture
Notes in Math., 1255, Springer, Berlin, 1987

\noindent
[BrS] $\underline{\; \;}$ ; Salamon, Simon M., On the construction of some complete
metrics with exceptional holonomy. Duke Math. J. 58 (1989), no. 3, 829-850

\noindent
[ChW] Chern, S. S. ; Wolfson, Jon G., Minimal surfaces by moving frames.
Amer. J. Math. 105 (1983), no. 1, 59--83

\vspace{3pc}

\noindent
Sung Ho Wang \\
C.R.M., Universite de Montreal \\
Montreal, Quebec H3C 3J7 \\
wang@crm.umontreal.ca

\end{document}